\documentclass[112pt]{amsart}
\usepackage{mathrsfs}
\usepackage{amsfonts}
\usepackage{amssymb}
\usepackage[papersize={6.7in,9.4in},textwidth=13.4cm,textheight=20cm,centering]{geometry}
\usepackage{enumerate}

    \usepackage[utf8]{inputenc}
    \usepackage[T1]{fontenc}
    \usepackage{libertine} 
    \usepackage[libertine]{newtxmath}

\usepackage{amsmath}
\numberwithin{equation}{section}

\usepackage{multirow}

\usepackage[colorlinks=true,citecolor=blue,linkcolor=blue]{hyperref}
\hypersetup{
pdfstartpage=1,
pdfstartview=FitH}

\usepackage{amsthm,amssymb}
\usepackage{indentfirst}
\usepackage[leqno]{amsmath}

\allowdisplaybreaks

\newtheorem*{theorem*}{Theorem}
\newtheorem{theorem}{Theorem}
\newtheorem{lemma}{Lemma}

\newtheorem*{claim*}{Claim}

\theoremstyle{definition}

\DeclareMathOperator{\Mod}{mod}

\renewcommand{\bmod}[1]{\,(\Mod{ #1})}

\newcommand{\R}{\mathbb{R}}

\newcommand{\balpha}{\boldsymbol{\alpha}}
\newcommand{\bbeta}{\boldsymbol{\beta}}
\newcommand{\bgamma}{\boldsymbol{\gamma}}

\newcommand{\bZ}{\mathbf{Z}}

\newcommand{\cL}{\mathcal{L}}

\newcommand{\cN}{\mathcal{N}}

\def\le{\leqslant}

\def\ge{\geqslant}

\begin{document}

\title{A problem of D. H. Lehmer in short intervals. II}
\author{Qixiang Shen}

\address{School of Mathematics and Statistics, Xi'an Jiaotong University, Xi'an 710049, P. R. China}
\email{qxshen@yeah.net}

\begin{abstract}
     A problem of D. H. Lehmer suggests to study the number of integers, each of which has different parity from its multiplicative inverse modulo $q.$
     For large prime $q$, we obtain an asymptotic formula for the number of such integers up to $N$, where $N$ is a bit smaller than $q^{1/2}$. This beats the barrier $q^{1/2}$ in the prime modulus case. An estimate for the second moment of the error term on average over $q$ is also established. The main inputs are estimates for several bilinear forms with Kloosterman fractions.
\end{abstract}

\maketitle


\section{Introduction and our main result}
Let $q>1$ be an odd integer. For any $a\in\bZ$ with $(a, q)=1$, let $\bar{a}$ denote the unique integer satisfying $1\le \bar{a}\le q-1$ and $a\bar{a}\equiv 1\bmod q$. Define
\begin{align*}
    \cL(q)=\{1\le a\le q-1:(a, q)=1, 2\nmid a+\bar{a}\}.
    \end{align*}
The problem of determining the size of $\mathcal{L}(q)$ was first proposed by D. H. Lehmer and recorded as Problem F12 in Guy’s book \cite{Gu81}. This problem was initated by Zhang \cite{Zha92,Zha93,Zha94}, who proved an asymptotic formula for $\vert\mathcal{L}(q)\vert$ when $q$ is a large odd integer.
We refer the readers to \cite{LRS07,LY09,Sh09} for related progress.

Now we continue our study on the short-interval version of the original problem of D. H. Lehmer. Let
\begin{align*}
    \cL(N,q)=[1,N]\cap \cL(q),\ \ \ L(N, q)=\vert\cL(N ,q)\vert.
\end{align*}
It is natural to expect, for all $N>q^\varepsilon$ with any $\varepsilon>0,$ that
\begin{align}\label{eq:goal}
    L(N, q)\sim \frac{1}{2}N\varphi(q)q^{-1}
\end{align}
holds with an implicit error term. Using Weil's bound for individual Kloosterman sums, Zheng \cite{Zhe93} proved that
\begin{align}\label{eq:zheng}
    L(N,q)=\frac{1}{2}N\varphi(q)q^{-1}+O(q^{1/2}\tau(q)\log^2q).
\end{align}
Note that (\ref{eq:zheng}) is nontrivial for all $N\gg q^{1/2+\varepsilon}$. In particular, it gives nothing nontrivial as long as $N\sim q^{1/2}.$ In a previous paper \cite{Sh26}, we proved \eqref{eq:goal} when $N\gg q^{\varepsilon}$ for any given $\varepsilon>0$, provided that $q$ is squarefree and has only prime factors not exceeding $q^{\eta}$, where $\eta>0$ is sufficiently small in terms of $\varepsilon$. Clearly, this result does not apply to prime moduli, in which case the $1/2$-barrier remains out of reach. 

In this paper, we first revisit this problem when $q$ is a large prime and $N$ is a bit smaller than $q^{1/2}.$

\begin{theorem}\label{thm:main}
Let $q$ be an odd prime, $1\le N\le q-1$. For any $0<\delta<1/4$, whenever \[N\ge q^{1/2}\exp\big(-(\log q)^{1/2-2\delta}\big),\] we have
    \begin{align}\label{eq:main result}
        L(N ,q)=\frac{1}{2}N\varphi(q)q^{-1}+O\Big(\frac{N}{(\log q)^{\delta}}\Big).
    \end{align}
The implied constant depends at most on $\delta$.
\end{theorem}

The proof of Theorem \ref{thm:main} relies on a deep estimate of Bourgain and Garaev \cite{BG14} on incomplete Kloosterman sums (Lemma \ref{lm:bourgain} below) and the "shift by $ab$" trick of Vinogradov, Burgess, Karatsuba and Friedlander--Iwaniec. More precisely, the latter one refers to the following double sum
\begin{align*}
\Sigma(\balpha)=\Sigma_k(\balpha;R,N,q)=\sum_{r\sim R}\alpha_r\sum_{n\le N}e\Big(\frac{kr\overline{n}}{q}\Big),
\end{align*}
where $\balpha=(\alpha_r)$ is a sequence of arbitrary coefficients. We also refer the readers to \cite{Sh19,KSWX23} for related progress.

\begin{theorem}\label{thm:kl sum estimate 1}
Let $q$ be a prime, $(k,q)=1$, $1<R<q/2$ and $R<N<q$. For any $s\in \mathbb{Z}^+$, we have
    \begin{align*}
    \Sigma(\balpha)\ll \Vert\balpha\Vert_{\infty}RN\log q\Big(\frac{s!q\log q}{NR(N/R)^s}+\frac{sq^{1/2}\log q}{NR}+\frac{s!}{(N/R)^s}+\frac{s}{q^{1/2}}\Big)^{1/2s}.
\end{align*}
The implied constant is absolute.
\end{theorem}

\begin{theorem}\label{thm:kl sum estimate}
Let $q$ be a prime, $(k,q)=1$, $1<R<q/2$ and $1<N<q$. For any $s\in \mathbb{Z}^+$, we have
    \begin{align*}
    \Sigma(\balpha)\ll  \Vert\balpha\Vert_{\infty}RN\log q\Big(\frac{sq^{1+1/s}\log q}{RN^2}+\frac{s}{N}\Big)^{1/2s}.
\end{align*}
The implied constant is absolute.
\end{theorem}

As mentioned above, we expect the asymptotic behaviour \eqref{eq:goal} is valid for all $N>q^\varepsilon$ with any $\varepsilon>0.$ However, the proof seems beyond the current approach for a general modulus $q$. We are able to show this is indeed the case on average over $q$. In fact, for a generic sequence $\balpha=(\alpha_n)$, we can prove that the difference
\begin{align*}
    \Delta(\balpha,N;q):=\sum_{n\in \cL(N, q)}\alpha_n-\frac{1}{2}\sum_{\substack{n\le N\\(n, q)=1}}\alpha_n
\end{align*}
is also small on average over $q$.

\begin{theorem}\label{thm:secondmoment}
Let $1<N<Q$ and $\balpha=(\alpha_n)$ an arbitrary sequence of coefficients. For any $\varepsilon>0$, we have
    \begin{align*}
        \sum_{\substack{q\sim Q\\ q\ {\rm odd}}}\vert \Delta(\balpha,N;q) \vert^2\ll\Vert\balpha\Vert_2^2 N^{11/12}Q^{1+\varepsilon}.
\end{align*}
The implied constant depends at most on $\varepsilon$.
\end{theorem}

The proof of Theorem \ref{thm:secondmoment} relies on estimates for a different bilinear form with Kloosterman fractions, due to Duke-Friedlander-Iwaniec \cite{DFI97} and Bettin-Chandee \cite{BC18}. See Lemmas \ref{lm:DFI} and \ref{lm:BC} below for details. Taking $\balpha\equiv 1$, Theorem \ref{thm:secondmoment} implies that, for any fixed $\eta>0$, if $N\gg Q^{\eta}$, then by choosing $0<\varepsilon<\eta/12$, (\ref{eq:goal}) holds for all but $o(Q)$ odd moduli $q\sim Q$.

\section{Preliminary lemmas}

The following lemma transforms the original problem to exponential sums. We refer the readers to \cite[Lemma 2]{Sh26} for the proof.

\begin{lemma}\label{lm:sum over symmetric set}
Let $q>1$ be an odd integer, $\{{\alpha_n}\}$ an arbitrary coefficient. Then we have
    \begin{align}\label{eq:transform into character sums}
        \sum_{n\in \cL (N, q)}\alpha_n=\frac{1}{2}\sum_{\substack{n\le N\\(n, q)=1}}\alpha_n+O\Big(\sum_{j=1}^2\sum_{0<\vert r\vert<q/2}\frac{1}{|r|}\Big\vert\sum_{n\in I_j}\alpha_{n,j}e\Big(\frac{ k_jr\overline{n}}{q}\Big)\Big\vert\Big),
        \end{align}
where $I_1=[1, N/2]$, $\alpha_{n, 1}=\alpha_{2n}$, $4k_1\equiv \pm1\bmod q$ and $I_2=[1, N]$, $\alpha_{n, 2}=\alpha_{n}$, $2k_2\equiv \pm1\bmod q$. The implied $O$-constant is absolute.
\end{lemma}

We need the following estimate for very short Kloosterman sums, which was established by Bourgain and Garaev \cite[Theorem 5]{BG14}. We refer the readers to \cite{Ko00, KR21} for related progress and to \cite{Ko10} for the weighted variants.

\begin{lemma}\label{lm:bourgain}
    Let $q>1$. For any fixed $c>0$, if $N>q^{c}$, then
    \begin{align*}
        \underset{(a, q)=1}{\max}\Big\vert \sum_{\substack{n\le N\\(n,q)=1}}e\Big(\frac{a\bar{n}}{q}\Big)\Big\vert<N\frac{(\log\log q)^{O(1)}}{(\log q)^{1/2}},
    \end{align*}
    where the implied constant depends at most on $c$. 
\end{lemma}

We recall the following bound on the number of solutions of multiplicative congruences, which was first studied by Ayyad, Cochrane and Zheng \cite[Theorem 1]{ACZ96}. We shall use the following stronger result of Kerr \cite[Theorem 1]{Ke17}.

\begin{lemma}\label{lm:CZ}
Let $q$ be a prime and $I=\{I_i\}_{i=1}^4$ a sequence of intervals not containing $0$ mod $q$. Define
\begin{align*}
   N(I)=\# \{(x_1, x_2, x_3, x_4)\in \prod_{i=1}^4I_i:  x_1x_2\equiv x_3x_4\bmod q\}.
\end{align*}
We have
\begin{align*}
        N(I)=\frac{\prod_{i=1}^4\vert I_i\vert}{q}+O(\prod_{i=1}^4\vert I_i\vert^{1/2}\log q).
\end{align*}
\end{lemma}

Thanks to Weil \cite{We48}, one has very strong bounds for one-dimensional exponential and character sums of rational function over finite fields, with squareroot cancellations. 
We shall use the following version, as a special case of \cite[Lemma 3.1]{MSW24}.
\begin{lemma}\label{lm:Weil bound}
    Let $q$ be a prime and $\mathbb{F}_{q}$ a finite field. Let $\psi$ be a non-principal additive character and let $g(X)\in \mathbb{F}_{q}(X)$ be a rational function of degree $d$ over $\mathbb{F}_{q}$. Assume that $g(X)\neq \alpha(h(X)^q-h(X))$ for all rational functions $h(X)\in\overline{\mathbb{F}_q}(X)$ and $\alpha\in \overline{\mathbb{F}_q}$. Then we have
    \begin{align*}
    \Big\vert\sum_{x\in\mathbb{F}_{q}}\psi(g(x))\Big\vert\le 2dq^{1/2},
    \end{align*}
where the sum is taken over those $x\in \mathbb{F}_q$ for which $g(x)$ is defined.
\end{lemma}

To prove Theorem \ref{thm:secondmoment}, we need a different bilinear form involving Kloosterman fractions. To be precise, let
\begin{align}\label{eq:bilinear form of kl sum}
    \mathscr{B}(\balpha, \bbeta)=\mathscr{B}_a(\balpha, \bbeta; M, N)=\sum_{m\sim M}\sum_{n\sim N}\alpha_m\beta_ne\Big(\frac{a\overline{m}}{n}\Big),
\end{align}
where $\balpha=(\alpha_m)$ and $\bbeta=(\beta_n)$ are two arbitrary sequences, and we keep $(m,n)=1$ in mind. The heart of the problem is to beat the trivial bound $\vert\mathscr{B}(\balpha, \bbeta)\vert\le\Vert\balpha\Vert_2\Vert\bbeta\Vert_2\sqrt{MN}$ under mild assumptions on $M,N.$
To this end, one can first apply Cauchy's inequality
\begin{align*}
    \vert\mathscr{B}(\balpha, \bbeta)\vert^2
    \le\Vert\balpha\Vert_2^2\sum_{m\sim M}\Big\vert\sum_{n\sim N}\beta_n e\Big(\frac{a\bar{m}}{n}\Big)\Big\vert^2
    =\Vert\balpha\Vert_2^2\mathop{\sum\sum}_{n_1,n_2\sim N}\beta_{n_1} \overline{\beta_{n_2} }\sum_{m\sim M}e\Big(\frac{a(n_2-n_1)\bar{m}}{n_1n_2}\Big).
\end{align*}
A naive approach is to use the P\'olya--Vinogradov method in the last $m$-sum, as in the proof of \cite[Theorem 1]{DFI97}. One then arrives at bounding complete Kloosterman sums following Weil, and this produces a nontrivial bound for $\vert\mathscr{B}(\balpha, \bbeta)\vert$ whenever $M\gg N^{1+\varepsilon}$. 

To obtain a nontrivial estimate in the balanced range $M\asymp N$, Duke, Friedlander and Iwaniec \cite[Theorem 2]{DFI97} adapted the amplification method developed in their earlier works \cite{DFI93, DFI94, DFI95}. As a result, in the case $a\ll MN$ and $M\asymp N$, they obtained a power-saving against the trivial bound. This idea was recently enhanced and refined by Bettin and Chandee \cite[Theorem 1]{BC18}. More generally, Bettin--Chandee's result allows an additional average over the parameter $a$. For our purpose in this note, we only present a special case of their result in the following as stated in Lemma \ref{lm:BC}.

\begin{lemma}[{\cite[Theorem 1]{DFI97}}]\label{lm:DFI}
For any integer $a\neq 0$ and any $\varepsilon>0$, we have
\begin{align}\label{eq:DFI}
    \mathscr{B}(\balpha, \bbeta)\ll\Vert\balpha\Vert_2\Vert\bbeta\Vert_2\Big((M+N)^{1/2}+\Big(1+\frac{\vert a\vert }{MN}\Big)^{1/2}\min(M, N)\Big)(MN)^{\varepsilon}.
\end{align}
    The implied constant depends at most on $\varepsilon$.
\end{lemma}

\begin{lemma}[{\cite[Theorem 1]{BC18}}]\label{lm:BC}
    For any integer $a\neq 0$ and any $\varepsilon>0$, we have
    \begin{align}\label{eq:BC}
 \mathscr{B}(\balpha, \bbeta)\ll\Vert\balpha\Vert_2\Vert\bbeta\Vert_2\Big(1+\frac{\vert a\vert }{MN}\Big)^{1/2}\Big((MN)^{7/20}(M+N)^{1/4}+(MN)^{3/8}(N+M)^{1/8}\Big)(MN)^{\varepsilon}.
    \end{align}
The implied constant depends at most on $\varepsilon$.
\end{lemma}

Combining Lemma \ref{lm:DFI} and Lemma \ref{lm:BC}, we obtain a non-trivial estimate whenever $M^{\eta}<N<M^{1/\eta}$, for any fixed $\eta>0$. More precisely, applying (\ref{eq:DFI}) if $MN\ll (M+N)^{25/13}$ and (\ref{eq:BC}) otherwise, we conclude the following estimate.

\begin{lemma}\label{lm:combine bound kl fraciton}
For any integer $a\neq 0$ and any $\varepsilon>0$, we have
\begin{align*}
    \mathscr{B}(\balpha, \bbeta)\ll\Vert\balpha\Vert_2\Vert\bbeta\Vert_2(\vert a\vert+ MN)^{1/2}(M+N)^{1/24}(MN)^{-1/24+\varepsilon}.
\end{align*}
The implied constant depends at most on $\varepsilon$.
\end{lemma}

The following duality principle allows us to understand $L^2$ estimate in terms of bilinear forms. This can be found in many references, and see \cite[Lemma 2]{Mo78} for instance.
\begin{lemma}\label{lm:duality principle}
The following three statements are equivalent:

$(a)$ For any complex numbers $\alpha_m,$
\[\sum_{n}\Big|\sum_m\alpha_m\phi(m,n)\Big|^2\leqslant\Delta\|\balpha\|_2^2.\]

$(b)$ For any complex numbers $\beta_n,$
\[\sum_{m}\Big|\sum_n\beta_n\phi(m,n)\Big|^2\leqslant\Delta\|\bbeta\|_2^2.\]

$(c)$ For any complex numbers $\alpha_m,$ $\beta_n,$
\[\Big|\sum_{m}\sum_n\alpha_m\beta_n\phi(m,n)\Big|\leqslant\sqrt{\Delta}\|\balpha\|_2\|\bbeta\|_2.\]
Here $\Delta$ is the same in all inequalities.
\end{lemma}

\section{Proof of Theorem \ref{thm:kl sum estimate 1}}
We introduce the shift $n\mapsto n+ur$, in the spirit of the "shift by $ab$" trick as in \cite{FI93, FM98, KMS17}, which implies that for some $\xi\in \R$, we have
\begin{align}\label{eq:after shift ab}
    \Sigma(\balpha)\ll \frac{\Vert\balpha\Vert_{\infty}\log q}{U}\sum_{r\sim R}\sum_{n\in \cN}\Big\vert \sum_{u\sim U} e\Big(\frac{kr(\overline{n+ur})}{q}\Big)e(\xi u)\Big\vert,
\end{align}
where $\cN$ is another interval of length at most $2N$, and $U$ is to be optimized later subject to the constraint $UR\le N$ and terms with $n+ur\equiv 0\bmod q$ are omitted.

To group variables, we put
\begin{align*}
    \nu(x)=\#\{(n, r): n\equiv xr\bmod q, n\in \cN, r\sim R\},
\end{align*}
for $x\,\mathrm{mod}\,q$. Thus we can write
\begin{align*}
    \Sigma(\balpha)\ll\frac{\Vert\balpha\Vert_{\infty}\log q}{U}\sum_{x\bmod q}\nu(x)\Big\vert \sum_{u\sim U} e\Big(\frac{k(\overline{x+u})}{q}\Big)e(\xi u)\Big\vert.
\end{align*}
By Hölder's inequality, we have for any $s\in \mathbb{Z}^+$, 
\begin{align}\label{eq:Holder 1}
    \Sigma(\balpha)\ll \frac{\Vert\balpha\Vert_{\infty}\log q}{U}\Sigma_1^{1-1/s}(\Sigma_2\Sigma_3)^{1/2s},
\end{align}
where
\begin{align*}
    &\Sigma_1=\sum_{x\bmod q} \nu(x),\\
    &\Sigma_2=\sum_{x\bmod q}\nu(x)^2,\\
    &\Sigma_3=\sum_{x\bmod q}\Big\vert \sum_{u\sim U} e\Big(\frac{k(\overline{x+u})}{q}\Big)e(\xi u)\Big\vert^{2s}.
\end{align*}

Clearly,
\begin{align}\label{eq:bound Sigma_1}
\Sigma_1\ll  RN.
\end{align}
And $\Sigma_2$ counts the number of solutions of
\begin{align*}
    n_1r_2\equiv n_2r_1\bmod q, n_i\in\mathcal{N}, r_i\sim R.
\end{align*}
By Lemma \ref{lm:CZ}, we have
\begin{align}\label{eq:bound Sigma_2}
    \Sigma_2\ll RN\Big(\log q+\frac{RN}{q}\Big).
\end{align}
If necessary, we split $\mathcal{N}$ into $O(1)$ intervals of length $<q$. The possible integers divisible by $q$ contribute $O(R^2)$, which is absorbed by $RN\log q$.

For $\Sigma_3$, we have
\begin{align*}
    \Sigma_3\le\sum_{u_1,\cdots, u_{2s}\le U}\Big\vert\sum_{x\bmod q}e\Big(\frac{k\big(\sum_{i=1}^s(\overline{x+u_i})-\sum_{i=s+1}^{2s} (\overline{x+u_{i}})\big)}{q}\Big)\Big\vert.
\end{align*}
If $\{u_1,\cdots,u_{s}\}=\{u_{s+1},\cdots,u_{2s}\}$, then the sum over $x$ is bounded by $q$ and thus the total contribution from such solutions is $O(s!U^sq)$. For other choices of $u_1,\cdots,u_{2s}$, we apply Lemma \ref{lm:Weil bound} with 
\begin{align*}
    g(X)=k\Big(\sum_{i=1}^s(X+u_i)^{-1}-\sum_{i=s+1}^{2s}(X+u_i)^{-1}\Big)=k\sum_a \frac{c_a}{X+a},
\end{align*}
where 
\begin{align*}
    c_a=\#\{1\le i\le s: u_i=a\}-\#\{s<i\le 2s: u_i=a\},
\end{align*}
and not all $c_a$ vanish. Hence, $g$ has at least one simple pole and at most $2s$ simple poles. In particular, $g$ cannot be of the form $\alpha(h(X)^q-h(X))$, since every pole of $h^q-h$, if any exists, has order divisible by $q$. Now Lemma \ref{lm:Weil bound} gives
\begin{align*}
    \left\vert\sum_{x\bmod q} e\Big(\frac{g(x)}{q}\Big)\right\vert\ll s q^{1/2},
\end{align*}
and the total contribution from such solutions is $O(sU^{2s}q^{1/2})$. Hence, 
\begin{align}\label{eq:bound Sigma_3}
    \Sigma_3\ll s!U^sq+sU^{2s}q^{1/2}.
\end{align}
Combining (\ref{eq:Holder 1}), (\ref{eq:bound Sigma_1}), (\ref{eq:bound Sigma_2}) and (\ref{eq:bound Sigma_3}) and choosing $U\asymp N/R$, we have
\begin{align*}
    \Sigma(\balpha)&\ll \frac{\Vert\balpha\Vert_{\infty}\log q}{U}(RN)^{1-1/2s}\Big(\log q+\frac{RN}{q}\Big)^{1/2s}\Big(s!U^sq+sU^{2s}q^{1/2}\Big)^{1/2s}\\
    &\ll \Vert\balpha\Vert_{\infty}RN\log q\Big(\frac{s!q\log q}{NR(N/R)^s}+\frac{sq^{1/2}\log q}{NR}+\frac{s!}{(N/R)^s}+\frac{s}{q^{1/2}}\Big)^{1/2s}.
\end{align*}
Then Theorem \ref{thm:kl sum estimate 1} follows.

\section{Proof of Theorem \ref{thm:kl sum estimate}}

To prove Theorem \ref{thm:kl sum estimate}, we introduce a different version of the "shift by $ab$" trick. For some $\xi\in \R$, we have
\begin{align*}
    \Sigma(\balpha)\ll \frac{\Vert\balpha\Vert_{\infty}\log q}{UV}\sum_{r\sim R}\sum_{n\in \cN}\sum_{v\sim V}\Big\vert \sum_{u\sim U} e\Big(\frac{kr(\overline{n+uv})}{q}\Big)e(\xi u)\Big\vert,
\end{align*}
where $\cN$ is another interval of length at most $2N$, and $U, V$ are to be optimized later subject to the constraint $UV\le N$.

To group variables, we put
\begin{align*}
    \gamma(x, y)=\#\{(n, r, v): kr\equiv xv\bmod q, n\equiv yv\bmod q, n\in \cN, r\sim R, v\sim V\},
\end{align*}
for $x, y\,\mathrm{mod}\,q$. Thus we can write
\begin{align*}
    \Sigma(\balpha)\ll\frac{\Vert\balpha\Vert_{\infty}\log q}{UV}\sum_{x, y\bmod q}\gamma(x, y)\Big\vert \sum_{u\sim U} e\Big(\frac{x(\overline{y+u})}{q}\Big)e(\xi u)\Big\vert,
\end{align*}

By Hölder's inequality, we have for any $s\in \mathbb{Z}^+$, 
\begin{align}\label{eq:Holder}
    \Sigma(\balpha)\ll \frac{\Vert\balpha\Vert_{\infty}\log q}{UV}W_1^{1-1/s}(W_2W_3)^{1/2s},
\end{align}
where
\begin{align*}
    &W_1=\sum_{x, y\bmod q} \gamma(x, y),\\
    &W_2=\sum_{x, y\bmod q}\gamma(x, y)^2,\\
    &W_3=\sum_{x, y\bmod q}\Big\vert \sum_{u\sim U} e\Big(\frac{x(\overline{y+u})}{q}\Big)e(\xi u)\Big\vert^{2s}.
\end{align*}

Clearly,
\begin{align}\label{eq:bound S_1}
W_1\ll  RNV.
\end{align}
And $W_2$ counts all solutions to the system of congruences
\begin{align*}
    &r_1v_2\equiv r_2v_1\bmod q,\\
    &n_1v_2\equiv n_2v_1\bmod q.
\end{align*}
in $r_1, r_2\sim R, n_1, n_2\in\mathcal{N}, v_1, v_2\sim V$. If necessary, we split $\mathcal{N}$ into $O(1)$ intervals of length $<q$. Once $(r_1, r_2, v_1, v_2, n_1)$ is fixed, the second congruence $n_1v_2\equiv n_2v_1\bmod q$ determines $n_2$ uniquely modulo $q$ and there is at most one admissible $n_2$. Hence, by Lemma \ref{lm:CZ}, we have
\begin{align}\label{eq:bound S_2}
    W_2\ll N\Big(\frac{(RV)^2}{q}+RV\log q\Big).
\end{align}

Applying the orthogonality of additive characters, we have
\begin{align*}
    &W_3\le q\sum_{u_1,\cdots, u_{2s}\le U}\sum_{\substack{y\nequiv -u_i\bmod q\\\ \sum_{i=1}^s \overline{y+u_i}\equiv\sum_{i=s+1}^{2s} \overline{y+u_{i}}\bmod q}}1.
\end{align*}
If $\{u_1,\cdots,u_{s}\}=\{u_{s+1},\cdots,u_{2s}\}$, then there are at most $q$ values for $y$ and thus the total contribution from such solutions is $O(s!U^sq)$. For other choices of $u_1,\cdots,u_{2s}$, there are at most $2s-1$ values of $y$ and thus the total contribution from such solutions is $O(sU^{2s})$. Hence, $W_3$ can be bounded by
\begin{align*}
    &\ll s!q^2U^s+sqU^{2s}.
\end{align*}
Choosing $U=(s!q/s)^{1/s}$ and using the Stirling's approximation, we have
\begin{align}\label{eq:bound S_3}
    W_3\ll s^{2s}q^{3}.
\end{align}

Combining (\ref{eq:Holder}), (\ref{eq:bound S_1}), (\ref{eq:bound S_2}) and (\ref{eq:bound S_3}) and choosing $V=N/U=N(s!q/s)^{-1/s}$, we have
\begin{align*}
    \Sigma(\balpha)&\ll \Vert\balpha\Vert_{\infty}RN\log q\Big(\frac{sq\log q}{RNV}+\frac{s}{N}\Big)^{1/2s}
    \ll  \Vert\balpha\Vert_{\infty}RN\log q\Big(\frac{sq^{1+1/s}\log q}{RN^2}+\frac{s}{N}\Big)^{1/2s}.
\end{align*}
Then Theorem \ref{thm:kl sum estimate} follows.

\section{Proof of Theorem \ref{thm:main}}
Applying Lemma \ref{lm:sum over symmetric set} with $a_n\equiv 1$, we have
\begin{align*}
    \sum_{\substack{a\le N\\(a, q)=1}}1=\sum_{d\vert q}\mu(d)\Big[\frac{N}{d}\Big]=N\varphi(q)q^{-1}+O(\tau(q)),
\end{align*}
and
\begin{align*}
    L(N, q)=\frac{1}{2}N\varphi(q)q^{-1}+O\Big(S(N, q)+\tau(q)\Big),
\end{align*}
where
\begin{align*}
    S(N ,q)=\sum_{j=1}^2\sum_{0<\vert r\vert<q/2}\frac{1}{\vert r\vert}\Big\vert\sum_{n\in I_j}e\Big(\frac{k_jr\bar{n}}{q}\Big)\Big\vert,
\end{align*}
with $I_1=[1, N/2]$, $4k_1\equiv \pm1\bmod q$ and $I_2=[1, N]$, $2k_2\equiv \pm1\bmod q$. Thus, it's sufficient to estimate $S(N, q)$, since $\tau(q)=2$. 

For $0<\delta<1/4$, define
\begin{align*}
    A=4(\log q)^{1/2-2\delta}.
\end{align*}
By the hypothesis, we have $N\ge q^{1/2}\exp(-A/4)$ throughout the proof. By a dyadic decomposition, we have
\begin{align}\label{eq:dyadic}
    S(N ,q)\ll \vert S_1(N, q, e^A)\vert+\log q \underset{e^A\le R<q/2}\max \vert S_2(N, q, R)\vert,
\end{align}
where
\begin{align*}
    &S_1(N, q, e^A)=\sum_{j=1}^2\sum_{0<\vert r\vert<e^A}\frac{1}{\vert r\vert}\Big\vert\sum_{n\in I_j}e\Big(\frac{k_jr\bar{n}}{q}\Big)\Big\vert,\\
    &S_2(N, q, R)=\frac{1}{R}\sum_{j=1}^2\sum_{\vert r\vert\sim R}\beta_{r,q}\sum_{n\in I_j}e\Big(\frac{k_jr\bar{n}}{q}\Big),
\end{align*}
with $\vert\beta_{r, q}\vert\ll 1$.
 
For $S_1(N, q, e^A)$, we apply Lemma \ref{lm:bourgain} to the innermost sum to get
\begin{align}\label{eq:final bound 1}
    S_1(N, q, e^A)\ll\frac{N(\log\log q)^{O(1)}}{(\log q)^{1/2}}\sum_{0<\vert r\vert<e^A}\frac{1}{\vert r \vert}\ll \frac{N(\log\log q)^{O(1)}A}{(\log q)^{1/2}}\ll_{\delta}\frac{N}{(\log q)^{\delta}},
\end{align}
since $N\ge q^{1/2}\exp(-A/4)>q^{1/3}$ for sufficiently large $q$.

For $S_2(N, q, R)$, we split the range $\vert r\vert\sim R$ into $r\sim R$ and $-r\sim R$. The negative range can be treated in the same way after replacing $k_j$ by $-k_j$. When $e^A\le R<N/q^{1/4}$, by Theorem \ref{thm:kl sum estimate 1}, we have 
\begin{align*}
    \underset{e^A\le R<N/q^{1/4}}{\max}S_2(N, q, R)&\ll N\log q\Big(\frac{s!q\log q}{NR(N/R)^s}+\frac{sq^{1/2}\log q}{NR}+\frac{s!}{(N/R)^s}+\frac{s}{q^{1/2}}\Big)^{1/2s}\\
    &\ll N\log q\Big(\frac{s!q^{1-(s-1)/4}\log q}{N^2}+\frac{sq^{1/2}\log q}{Ne^A}+\frac{s!}{q^{s/4}}+\frac{s}{q^{1/2}}\Big)^{1/2s}.
\end{align*}
Taking $s=2$, we get
\begin{align}\label{eq:final bound 2}
    \underset{e^A\le R<N/q^{1/4}}{\max}S_2(N, q, R)\ll N\log q\Big(\exp\big(-\frac{\Delta-\log\log q}{4}\big)+q^{-1/20}\Big),
\end{align}
where
\begin{align*}
    \Delta=A-\log\frac{q^{1/2}}{N}.
\end{align*}
Since $N\ge q^{1/2}\exp(-A/4)$, we have $\Delta\ge 3A/4$. Hence,
\begin{align*}
    \underset{e^A\le R<N/q^{1/4}}{\max}S_2(N, q, R)\ll N(\log q)^{5/4}\exp(-3A/16)\ll\frac{N}{(\log q)^{2026}}.
\end{align*}
When $N/q^{1/4}<R< q/2$, by Theorem \ref{thm:kl sum estimate}, we have
\begin{align*}
    \underset{N/q^{1/4}\le R< q/2}{\max}S_2(N, q, R)\ll  N\log q\Big(\frac{s^{1/2s}q^{5/4+1/s}\log q}{N^3}+\frac{s}{N}\Big)^{1/2s}.
\end{align*}
Taking $s=20$, we get
\begin{align}\label{eq:final bound 3}
    \underset{N/q^{1/4}\le R< q/2}{\max}S_2(N, q, R)\ll  N\log q\Big(\frac{q^{13/10}\log q}{N^3}+\frac{1}{N}\Big)^{1/40}\ll\frac{N}{(\log q)^{2026}} .
\end{align}

Combining (\ref{eq:dyadic}), (\ref{eq:final bound 1}), (\ref{eq:final bound 2}) and (\ref{eq:final bound 3}), we have
\begin{align*}
    S(N, q)&\ll_{\delta} \frac{N}{(\log q)^{\delta}},
\end{align*}
Then Theorem \ref{thm:main} follows.

\section{Proof of Theorem \ref{thm:secondmoment}}
By Lemma \ref{lm:sum over symmetric set}, we have
\begin{align*}
     \Delta(\balpha,N;q)\ll \sum_{j=1}^2\sum_{0<\vert r\vert<q/2}\frac{1}{|r|}\Big\vert\sum_{n\in I_j}\alpha_{n,j}e\Big(\frac{ k_jr\overline{n}}{q}\Big)\Big\vert,
\end{align*}
where $I_1=[1, N/2]$, $\alpha_{n, 1}=\alpha_{2n}$, $4k_1\equiv \pm1\bmod q$ and $I_2=[1, N]$, $\alpha_{n, 2}=\alpha_{n}$, $2k_2\equiv \pm1\bmod q$.
By a dyadic argument, we have
\begin{align*}
         \Delta(\balpha,N;q)\ll\log Q \sum_{j=1}^2\sum_{0<\vert r\vert<Q}\frac{1}{|r|}\underset{M\ll N}{\max}\Big\vert\sum_{m\sim M}\gamma_{m,j}e\Big(\frac{\pm r\overline{m}}{q}\Big)\Big\vert,
\end{align*}
where
\begin{align*}
    \gamma_{m, j}=\begin{cases}
        \alpha_{m/d_j, j} & d_j\vert m, m/d_j\in I_j,\\
        0 & \text{otherwise},
    \end{cases}
\end{align*}
with $d_1=4$ and $d_2=2$. And it's clear that for each 
$j$, $\Vert\bgamma_{\cdot, j}\Vert_2\ll\Vert \balpha\Vert_2$.

By Minkowski's inequality, we have
\begin{align}\label{eq:L2 norm}
    \Big(\sum_{\substack{q\sim Q\\ q\ {\rm odd}}}\vert\Delta(\balpha,N;q)\vert^2\Big)^{1/2}\ll \log Q \sum_{j=1}^2\sum_{0<\vert r\vert<Q}\frac{1}{|r|}\underset{M\ll N}{\max} S(M, r,j)^{1/2},
\end{align}
where
\begin{align*}
    S(M, r,j)=\sum_{\substack{q\sim Q\\ q\ {\rm odd}}}\Big\vert\sum_{m\sim M}\gamma_{m,j}e\Big(\frac{\pm r\overline{m}}{q}\Big)\Big\vert^2.
\end{align*}
And by Lemma \ref{lm:duality principle}, it's enough to prove the upper bound for
\begin{align*}
    \widetilde{S}(M, r,j)=\sum_{q\sim Q}\beta_q\sum_{m\sim M}\gamma_{m,j}e\Big(\frac{\pm r\overline{m}}{q}\Big),
\end{align*}
for any complex sequence $\bbeta=(\beta_q)$ supported on odd integers $q\sim Q$. 
Applying Lemma \ref{lm:combine bound kl fraciton} with $\varepsilon\leftarrow \varepsilon/4$, we have
\begin{align*}
    \widetilde{S}(M, r, j)&\ll\Vert\bgamma_{\cdot, j}\Vert_2\Vert\bbeta\Vert_2(\vert r\vert+MQ)^{1/2}(M+Q)^{1/24+\varepsilon/4}(MQ)^{-1/24}\\
    &\ll\Vert\balpha\Vert_2\Vert\bbeta\Vert_2 M^{11/24}Q^{1/2+\varepsilon/4},
\end{align*}
since $r\ll Q$ and $M\ll Q$. Hence, by Lemma \ref{lm:duality principle},
\begin{align*}
   \underset{M\ll N}{\max} S(M, r,j)\ll  \Vert\balpha\Vert_2^2N^{11/12}Q^{1+\varepsilon/2}.
\end{align*}
Finally, by (\ref{eq:L2 norm}), we have
\begin{align*}
     \sum_{\substack{q\sim Q\\ q\ {\rm odd}}}\vert\Delta(\balpha,N;q)\vert^2\ll \log^4 Q\Vert\balpha\Vert_2^2N^{11/12}Q^{1+\varepsilon/2}\ll \Vert\balpha\Vert_2^2N^{11/12}Q^{1+\varepsilon}.
\end{align*}
We are done!

\section{Acknowledgments}
The author thanks his advisor, Professor Ping Xi, for introducing this problem and useful discussions. This work is supported in part by Shaanxi NSF (No. 2025JC-QYCX-002), Shaanxi Fundamental Science Research Project for Mathematics and Physics (No.25JSZ007) and China Scholarship Council.

\end{document}